\documentclass[twocolumn, 10pt,pre]{revtex4-1}

\usepackage{amsmath}
\usepackage{amssymb}
\usepackage{amsthm}
\usepackage{graphicx}
\usepackage{subfig}

\newcommand{\p}{\partial}

\newcommand{\cf}{{\it cf}.~}

\newcommand{\sech}{\mbox{sech}}

\begin{document}

\title{Propagation of fronts in the Fisher-Kolmogorov equation with spatially varying diffusion}

\author{Christopher W. Curtis}
\affiliation{Department of Applied Mathematics, University of Colorado, Boulder, Colorado 80309, USA}
\author{David M. Bortz}
\affiliation{Department of Applied Mathematics, University of Colorado, Boulder, Colorado 80309, USA}

\begin{abstract}
The propagation of fronts in the Fisher-Kolmogorov equation with spatially varying diffusion coefficients is studied.  Using coordinate changes, WKB approximations, and multiple scales analysis, we provide an analytic framework that describes propagation of the front up to the minimum of the diffusion coefficient.  We also present results showing the behavior of the front after it passes the minimum.  In each case, we show that standard traveling coordinate frames do not properly describe front propagation.  Lastly, we provide numerical simulations to support our analysis and to show, that around the minimum, the motion of the front is arrested on asymptotically significant timescales. 
\end{abstract}
\maketitle
\section{Introduction}
Front propagation in reaction-diffusion equations (RDEs) is an important topic in the physical and biological sciences.  These models appear throughout biology \cite{murray,Murray2003}, ecology \cite{shigesada}, cancer research \cite{Maini2004}, chemical kinetics \cite{epstein}, and geochemistry \cite{grindrod}.  While classically concerned with homogeneous environments, interest in the behavior of fronts in heterogeneous environments has increased over the last several years.  In particular, traveling waves, or propagating fronts, through certain classes of time \cite{HammondBortz2011} and spatially varying environments \cite{shigesada,engel} have been well studied.  Nonlinear, density-dependent diffusion, where the diffusion coefficient depends on $u$, has also been thoroughly examined in \cite{murray,Murray2003} and \cite{grindrod} amongst others.  The related problem of varying selectivity, or excitability, was examined in \cite{engel,cuesta,keener,mendez}.   For the reader interested in more mathematical issues concerning front propagation, we refer to \cite{xin} for a thorough review of this topic.  Note that a complete listing of all references addressing these topics is beyond the scope of this paper, but the aforementioned references provide extensive bibliographies to the broader literature on front propagation in RDEs.  

We consider the case of a heterogeneous Fisher-Kolmogorov (FK), or Kolmogorov-Petrovskii-Piskunov (KPP), equation where heterogeneities are represented by spatially varying diffusion, {\it i.e.},
\begin{equation}
u_{t} = \left(a(x) u_{x}\right)_{x} + f(u).\label{FKax}
\end{equation}
In the biological literature, spatially-heterogeneous diffusion coefficients $a(x)$ are discussed repeatedly, but continuously varying models are rarely investigated rigorously.  For example, in ecological applications, several authors acknowledge the significance of heterogeneous diffusion \cite{Andow1990,Cruywagen1996,Dobzhansky1979,Hastings1982,Turchin1993} and some have developed approaches for estimating $a(x)$ across different habitats \cite{Dobzhansky1979,Harrison1980}.  However, many researchers dismiss spatially varying diffusion by attributing the diffusion variability to an evolutionary response in a sub-population (\cf \cite{Slatkin1985,Bridle2001}). When researchers do allow diffusion to vary, it is typically assumed to vary periodically in a square wave (\cf \cite{Cruywagen1996,Kinezaki2003,Kinezaki2010,shigesada,Shigesada1986}) as would commonly be encountered in agricultural or urban contexts.  In the spatially discrete, or lattice, context, researchers consider patches where diffusion is constant on each patch (\cf \cite{Petrovskii2005} and the references therein).

A morhpogenesis phenomenon in which $a(x)$ has been given serious consideration is in a general RDE model for the slug stage of {\it Dictyostelium discoideum}.  In this case, diffusion is directly modulated by the spatial distribution of a morphogen gradient concentration.  The resulting diffusion profile across a gap junction is a hyperbolic cosine, {\it i.e.}, $a(x)=\cosh (x)$.  The RDE model with this form of $a(x)$ has been investigated from an analytical (establishing a scale-invariant property for the generated wave pattern) \cite{Othmer1980} as well as computational \cite{Philip1992,Benson1993} perspective.  More recent work has considered the impact of letting $a(x)=D+ \eta x^2$ (for $D$ constant and $\eta$ small) in the context of the avascular growth phase of cancer \cite{Chaplain1994} as well as on Turing bifurcations of standing wave solutions to RDEs \cite{Benson1998}.  

None of these investigations, however, have generated analytical results for the impact on traveling waves of a continuously, strongly varying, diffusion coefficient.  Thus the form of $a(x)$ (motivated by the biologically justified choice of $\cosh(x)$) we choose to study is 
\begin{equation}
a(x) =  x^{2}+\epsilon.
\label{choice}
\end{equation}
From a modeling perspective, the diffusion coefficient represents an environment in which the necessary resources sustaining diffusion of the front dwindle and then grow as the front moves from left to right.  Thus, it offers a good canonical model for a diffusion-mediated barrier.  For example, one could use this mathematical structure to model a geographic barrier such as a mountain.  Further, this choice also allows us to study the effect of strongly varying diffusion, which is in contrast to the case looked at in \cite{mendez}.  This is because the choice in \eqref{choice} is more analytically tractable than the choice of hyperbolic cosine while still maintaining the concave shape with a global minimum.  In this work, we show, the strong variation creates several different asymptotic regimes through which the behavior of the front changes significantly.  Lastly, it is in some sense an examination of the solution behavior on a subdomain of the periodically varying case studied in \cite{shigesada}.  

Throughout the paper, we take as initial conditions for \eqref{FKax} the step initial condition
\begin{equation}
u(x,0)=\begin{cases}
1, & \, x\leq x_{c}(0)\\
0, & \, x>x_{c}(0)
\end{cases}\, ,
\label{step}
\end{equation}
where $x_{c}(0)$ denotes the location of the step.  It is well known in the traveling wave literature (\cite{aron}) that for this initial condition (and constant diffusion), the solution converges to a propagating front with speed determined by whether the nonlinearity induces pushed or pulled dynamics.  Initial conditions with slower rates of decay generate a waveform with a faster wave speed.  We label the {\it location} of the front via $x_{c}(t)$ with the convention that $u(x_{c}(t),t)=1/2$, so that the initial step is given by $u=1$ for $x < x_{c}(0)$ and $u=0$ for $x > x_{c}(0)$.  We then look for fronts, {\it i.e.}, solutions $u(x,t)$ bounded between zero and one that propagate from left to right along the spatial coordinate.  

In the FK equation, the minimum of the diffusion coefficient represents a {\it turning point}, {\it i.e.}, the sign of $a'(x)$ changes through the minimum. We develop, via multiple scales and WKB methods, an asymptotic description of the propagation of the front up to the point it crosses the turning point.  Our analysis shows because of the choice of a strongly varying diffusion coefficient that the front does not propagate in a simple traveling coordinate system, {\it i.e.}, $u(x,t)$ cannot be written in the form $u(x+ct)$.  Instead, the front travels along curves of the form $xe^{ct}=\tilde{x}_{0}$, where $\tilde{x}_{0}$ is some constant.  This is to say that we can write $u(x,t) = u(\ln(x) + ct)$, and we can then think of the front as a traveling wave in a more generalized sense.   

As can be seen from the analysis in Section \ref{sec:SFC}, the minimum of \eqref{choice} forces a fundamental shift in the behavior of the tail of $u(x,t)$ as the tail crosses the turning point.  This ultimately impacts the front by causing a shock like structure to form as it approaches the turning point.  At the turning point, we show numerically that the propagation of the front is arrested on a timescale of $\mathcal{O}(1/\sqrt{\epsilon})$.  This reflects the fact that at the turning point there is an asymptotically small amount of resources necessary for front propagation, and thus the front slows or is {\it trapped}. 

Beyond the turning point, we have preliminary analytic results which explain the behavior of the front.  We hypothesize that the front past the turning point is described by solutions to a stationary FK equation, in a traveling coordinate of the form $\ln(x)+ct$.  This again shows the need for a more general notion of traveling waves in the case of strongly varying diffusion.  However, more work is necessary to create a complete analytic treatment of the front beyond the turning point, especially for how the speed of the front is selected.  This will be addressed in a later paper.  We show via numerical simulation that our analytic treatment of the front up to the turning point is accurate.  The simulations also show, as analytically predicted, that there is a fundamental shift in the dynamics on either side of the turning point, and the numerics shows how fronts develop past the turning point.  

The structure of the paper is as follows.  In Section \ref{softfront}, we present our argument for how the front propagates up to the turning point.  In Section \ref{Softening IC}, we present an argument which shows how the FK equation smooths out step initial data.  In Section \ref{Reduction to a Stationary Equation}, we present a means of finding a traveling coordinate in which the FK equation is stationary, and we then present a series of arguments which supports the idea that beyond the turning point, fronts evolve according to this stationary equation.  Finally, in Section \ref{Numerical Results}, we present our numerical results.

\section{\label{sec:SFC}The Soft Front Approximation}
\label{softfront}When $x_{c}(t)\ll -1$, we assume that the solution $u$ has a {\it soft} front, {\it i.e.}, a condition we define via the asymptotic relationship 
\[
a(x)u_{xx} \ll a'(x)u_{x}.
\]  
This asymptotic condition is denoted as the Soft Front Approximation (SFA).  When the location of the front is such that $|x_{c}(t)|\gg 1$, we have that $a(x) \sim x^{2}$.  Thus, if for $\lambda > 0$, $u \sim 1/x^{\lambda}$ in a region around the front, then 
\[
a(x)u_{xx} \sim \lambda(\lambda+1)x^{-\lambda}, ~ a'(x)u_{x} \sim 2\lambda x^{-\lambda},
\]
so that the soft front condition requires that $0< \lambda \ll 1$.  We see that this approximation is not valid for step or rapidly decaying initial conditions.  This issue is addressed in the next section, in which we demonstrate how the front softens from step initial conditions.  As we will show, diffusion relaxes the rapidly decaying profile on short time scales and then makes the soft front approximation valid.  

Using the SFA, we now must solve the semi-linear hyperbolic equation
\[
u_{t} \sim  a'(x)u_{x} + f(u).
\]    
The method of characteristics then gives the system of differential equations
\[
\frac{dx}{dt}=  -a'(x), ~ \frac{du}{dt}  =  f(u).
\]
Using our choice of diffusion model $a(x)=x^{2}+\epsilon$, and choosing $f(u)=u(1-u)$, we get the solution 
\begin{equation}
u(x,t) \sim \frac{u(xe^{2(t-t_{0})},t_{0})e^{t-t_{0}}}{1 +u(xe^{2(t-t_{0})},t_{0})(e^{t-t_{0}}-1) },
\label{sfasol}
\end{equation}
which follows the characteristics $x(t) = x_{0}e^{-2(t-t_{0})}$.  In the case that $x_{0}<0$, we see that information propagates towards the origin as desired for a traveling front.  However, if $x_{0}>0$, then all information again propagates to the origin, and thus the soft front model cannot describe a front propagating past the origin.  Hence, we have a mechanism that explains front propagation that follows a decreasing diffusion coefficient, {\it i.e.}, $a(x)$ is strictly decreasing when $x_{c}(t)<0$.  

Therefore, given a snapshot of the profile $u(x,t)$ at a time $t_{0}$ where the SFA holds, we can describe the propagation of the front up to a neighborhood of the minimum of $a(x)$.  As can be seen from \eqref{sfasol}, the front does not propagate as a traveling wave, {\it i.e.}, there is not a constant wave speed such that $u(x,t) \sim u(x-ct)$.  However, if we suppose that 
\[
u(xe^{2(t-t_{0})},t_{0}) \sim 1, 
\]    
for $xe^{2(t-t_{0})}<x_{c}(t_{0})$, then we see that \eqref{sfasol} gives us that $u(x,t) \sim 1$.  Likewise, if we suppose that 
\[
u(xe^{2(t-t_{0})},t_{0}) \sim 0, 
\]    
for $xe^{2(t-t_{0})}>x_{c}(t_{0})$, then we see that \eqref{sfasol} gives us that $u(x,t) \sim 0$.  In this case then, taking logarithms, we get that the front travels along curves of the form 
\begin{equation}
\ln(x(t)) + 2(t-t_{0}) = \ln(x_{c}(t_{0})).
\label{twc}
\end{equation}
In this log-transformed spatial coordinate we can see the front propagates with speed $c = 2$.  We refer to \eqref{twc} as a {\it traveling wave coordinate} (TWC), and we see that the TWC and the SFA allow us to compute a generalized notion of wave speed.  This idea is studied further in Section \ref{sec:BeyondBarrier}, where we extend our analysis beyond the turning point.   

We note that on a finite time scale we do not imagine that $u(x,t)$ transitions globally to an algebraically decaying profile.  We are therefore arguing that the SFA holds in a region around $x_{c}(t)$ and this region must be matched, via intermediate layers, to the far field of $u(x,t)$ which should have a much steeper decay profile.  To begin to solve this problem, we look ahead of the front where $0<u\ll 1$.  To study this regime, we linearize the FK equation around $u=0$, {\it i.e.}, we let $u = \tilde{\epsilon} v$, and then collect all terms in $\tilde{\epsilon}$.  As is common, we assume that $f(u)>0$, $f(0)=f(1)=0$, $f'(0)=1$, and $f'(1)<0$, so that we get an equation for $v(x,t)$ of the form
\begin{equation}
v_{t} = \left(a(x)v_{x}\right)_{x} + v.
\label{linear}
\end{equation}
We suppose $v$ is given by the WKB ansatz, (\cf \cite{cuesta}, \cite{mendez}, and \cite{haberman} for examples of this approach)
\[
v(x,t) = A(x,t)e^{\phi(x,t)},
\]
from which we get the leading order problem
\begin{equation}
\phi_{t} + a(x)\phi^{2}_{x} + 1 = 0.  
\label{wkbphase}
\end{equation}
Using the method of characteristics and noting that \eqref{wkbphase} is a Hamilton-Jacobi equation (\cf \cite{evans}), with Hamiltonian $H(p,x)=1+a(x)p^{2}$, we get that 
\[
\frac{dx}{dt} = \pm 2 \tilde{H}\sqrt{a(x)}, ~ \phi(x,t) = (\tilde{H}^{2}-1)t + \phi_{0}(x_{0}(x,t)),
\]
where $\tilde{H}=\sqrt{H-1}$ is a constant along the characteristics.  The characteristics associated with the choice of diffusion coefficient \eqref{choice} are given by
\[
x_{0}+\sqrt{x_{0}^{2}+\epsilon} = \left(x+\sqrt{x^{2}+\epsilon} \right)e^{\mp 2\tilde{H}t}.
\] 
This expression is equivalent to 
\[
-\epsilon = \left(x_{0}-\sqrt{x_{0}^{2}+\epsilon}\right)\left(x+\sqrt{x^{2}+\epsilon} \right)e^{\mp 2\tilde{H}t},
\]
and we therefore get that 
\[
x_{0} = \frac{1}{2}\left(\left(x+\sqrt{x^{2}+\epsilon} \right)e^{\mp 2\tilde{H}t} - \frac{\epsilon}{x+\sqrt{x^{2}+\epsilon} }e^{\pm 2 \tilde{H}t}\right).
\]
Thus we see that for $|x|\gg \sqrt{\epsilon}$, the characteristics are to leading order given by $x_{0}\sim xe^{\mp2\tilde{H}t}$.  In general we see that information can propagate to or away from the origin.  Therefore, ahead of the front, the decaying tail of $u(x,t)$ propagates along characteristics of the form $x\sim x_{0}e^{-2\tilde{H}t}$ for that portion of the tail between the front and the turning point.  The portion of the tail beyond the turning point follows characteristic curves of the form $x\sim x_{0}e^{2\tilde{H}t}$ so that information travels away from the origin.  Note, the value of $\tilde{H}$ changes along different characteristics, and so we see that information propagates at different speeds along different characteristic curves.

Through the turning point, or when $|x| \ll \sqrt{\epsilon}$, we get the leading order behavior 
\[
x_{0} \sim \mp\sqrt{\epsilon}\sinh(2\tilde{H}t)+ \left(x+\frac{x^2}{2\sqrt{\epsilon}} \right)\cosh(2\tilde{H}t), 
\]
so that for $x_{0} < 0$ we have 
\[
x(t) \sim \sqrt{\epsilon}\left(-1 + \sqrt{1-2\tanh(2\tilde{H}t)+2\frac{x_{0}}{\sqrt{\epsilon}}\sech(2\tilde{H}t)}\right), 
\]
and for $x_{0}>0$ we have 
\[
x(t) \sim \sqrt{\epsilon}\left(-1 + \sqrt{1+2\tanh(2\tilde{H}t)+2\frac{x_{0}}{\sqrt{\epsilon}}\sech(2\tilde{H}t)}\right).
\]
Note, neither formula is useful on all time scales, but they are useful on times scales of $\mathcal{O}(1)$.  We see for $x_{0}<0$ that the characteristics in the layer $|x|\ll \sqrt{\epsilon}$ propagate away from the origin back to $-\sqrt{\epsilon}$.  Thus the characteristics in the outer region $|x|\gg \sqrt{\epsilon}$ collide in some intermediary layer with those in the inner region.  The WKB analysis thus shows how the SFA eventually breaks down.    

Likewise, the WKB analysis shows that as the front approaches the origin something like a shock, or steepened front, must form.  Determining this shock structure is beyond the scope of the current paper, though we are able to demonstrate it numerically in Section \ref{Numerical Results}.  We can get some hint though of the dynamics of the front from the following argument.  Once the front has entered the inner region, {\it i.e.}, $|x_{c}(t)|\ll \sqrt{\epsilon}$, the FK equation becomes  
\[
u_{t} \sim \epsilon u_{xx} + f(u).
\]
If the solution $u$ follows the dynamics of a propagating front, then $u$ is of the form
\[
u(x,t) \sim u\left(\frac{x-\sqrt{\epsilon}ct}{\sqrt{\epsilon}} \right),  
\]
so that we expect the front to be trapped in the inner region around the turning point on timescales of $\mathcal{O}(1/\sqrt{\epsilon})$.  While our analysis at this time is not complete, the numerical results supports this hypothesis.      

\section{\label{Softening IC}Softening Sharp Initial Conditions}
In this section we explain how starting with step initial conditions \eqref{step}, we can transition from the sharp to soft front regime.  We note that in the life of the traveling wave, the softening discussed here happens before the SFA becomes valid.  However, presenting it first would distract from the more important results in Section \ref{sec:SFC}.  We present it now so that a complete description of the behavior of the wave before barrier transit is available for reference in Section \ref{sec:BeyondBarrier}.

We assume throughout this section that $x_{c}(0)\gg1$.  The asymptotic condition describing the sharp front is that 
\[
a(x)u_{xx} \gg a'(x)u_{x},
\]
which, if $u\sim 1/x^{\lambda}$, implies that $\lambda \gg 1$.  We also choose a parameter $L\gg x_{c} \gg1$ so that $u_{x}(\pm L, t) = 0$.  We choose Neumann boundary conditions to allow for analytical tractability.  We introduce the fast time $T = \frac{t}{\epsilon}$, so that, using the ansatz,
\[
u = u_{0}(x,T,t) + \epsilon u_{1}(x,T,t) + \cdots,
\]  
we get the equations
\[
\begin{array}{rcl}
\p_{T}u_{0} & = &  \p_{x}\left(a(x)\p_{x}u_{0}\right), \\
&\\
\p_{T}u_{1} & = & \p_{x}\left(a(x)\p_{x}u_{1}\right) +f(u_{0}) - \p_{t}u_{0}.
\end{array}
\]
As for the leading order behavior $u_{0}$, we need only solve the linear diffusion equation to find it. Using separation of variables in space and time, we write $u_{0}(x,t,T)=\phi(x;t)\psi(T)$ which leads to the expansion for $u_{0}$ 
\begin{equation}
u_{0}(x,T,t)=\sum_{n=0}^{\infty}\sigma_{n}(t)\phi_{n}(x)e^{\lambda_{n}T},\label{loapprox}
\end{equation}
 where $\phi_{n}$ and $\lambda_{n}$ solve the Sturm-Liouville problem
\begin{equation}
 \p_{x}\left(a(x)\p_{x}\phi_{n}(x)\right)=\lambda_{n}\phi_{n},\;\p_{x}\phi_{n}(\pm L)=0.\label{evalproblem}
\end{equation}
 Sturm-Liouville theory (\textit{cf.} \cite{guenther}) ensures the functions $\phi_{n}$ form a complete, orthonormal set with respect to the norm 
\[
\left\Vert \phi\right\Vert ^{2}=\int_{-L}^{L}\left|\phi(x)\right|^{2}dx.
\]
 Since we assume that $a(x)>0$, then 
\begin{align*}
\lambda_{n}\int_{-L}^{L}\left|\phi_{n}(x)\right|^{2}dx& = &\int_{-L}^{L}\phi_{n}(x)\p_{x}\left(a(x)\p_{x}\phi_{n}\right)dx\\
&=&-\int_{-L}^{L}a(x)\left|\phi_{n,x}\right|^{2}dx\leq0,
\end{align*}
 so that $\lambda_{n}\leq0$. Further, we see $\phi_{0}$ is a constant corresponding to the eigenvalue $\lambda_{0}=0$. We set $\phi_{0}=\sqrt{1/2L}$ so that $||\phi_{0}||=1$. Likewise, using the step initial condition for $u$, we can find the initial condition for $\sigma_{0}(t)$ as 
\[
\sigma_{0}(0)=\frac{x_{c}+L}{\sqrt{2L}}.
\]
 The remaining terms $\sigma_{n}(t)$ have the initial conditions
\[
\sigma_{n}(0)=\frac{a(x_{c})}{\lambda_{n}}\p_{x}\phi_{n}(x_{c}/\epsilon).
\]

Moving to the second term $u_{1}(x,T,t)$, we see, using Duhamel's principle, that we can write $u_{1}$ as 
\[
u_{1}(x,T,t)=\int_{0}^{T}\sum_{n=0}^{\infty}\gamma_{n}(t,s)\phi_{n}(x)e^{-\left|\lambda_{n}\right|(T-s)}ds,
\]
 where 
\begin{equation}
\gamma_{n}(t,s)=\int_{-L}^{L}(-\partial_{t}u_{0}+u_{0}(1-u_{0}))\phi_{n}dx.\label{dncoeff}
\end{equation}
 Since $\lambda_{0}=0$, we see that a possibility for a secularity, which means the asymptotic series becomes invalid on $\mathcal{O}(1)$ time scales
(\textit{cf.} \cite{Bender}), arises from computing 
\[
\int_{0}^{T}\gamma_{0}(t,s)ds,
\]
 since if $\gamma_{0}(t,s)$ were independent of $s$ then $u_{2}$ would have a term that growing linearly in $T$. Expanding the integrand in Equation \eqref{dncoeff} gives 
\begin{align*}
-\partial_{t}u_{0}(x,t,s)+u_{0}(x,t,s)(1-u_{0}(x,t,s))= \\
 \phi_{0}(-\dot{\sigma}(t)+\sigma_{0}(t)(1-\phi_{0}\sigma_{0}(t)))\\
+ \sum_{n=1}^{\infty}((1-2\sigma_{0}\phi_{0})\sigma_{n}-\dot{\sigma}_{n})\phi_{n}e^{-\left|\lambda_{n}\right|s}\\
 -\sum_{n,j>0}^{\infty}\sigma_{n}\sigma_{j}\phi_{n}\phi_{j}e^{-(\left|\lambda_{n}\right|+\left|\lambda_{j}\right|)s}.
\end{align*}
 Using the orthonormality of the functions $\phi_{n}$, then from above one has 
\begin{align*}
\int_{0}^{T}\gamma_{0}(t,s)ds &= &T(-\dot{\sigma}_{0}(t)+\sigma_{0}(t)(1-\phi_{0}\sigma_{0}(t)))\\
&&+\phi_{0}\sum_{n=1}^{\infty}\frac{\sigma_{n}^{2}}{2\left|\lambda_{n}\right|}\left(e^{-2\left|\lambda_{n}\right|T}-1\right).
\end{align*}
 Thus, in order to remove the secularity, we enforce the condition
\[
\dot{\sigma}_{0}=\sigma_{0}(1-\phi_{0}\sigma_{0}),
\]
 which has the solution 
\[
\sigma_{0}(t)=\frac{1}{\phi_{0}+\tilde{\sigma}e^{-t}},
\]
 where $\tilde{\sigma}=-\phi_{0}+1/\sigma_{0}$.

As for the terms $\gamma_{n}$, one has 
\begin{align*}
\int_{0}^{T}\gamma_{n}(t,s)e^{-\left|\lambda_{n}\right|(T-s)}ds=  T((1-2\sigma_{0}\phi_{0})\sigma_{n}-\dot{\sigma}_{n})e^{-\left|\lambda_{n}\right|T}\\
 -\sum_{m,j>0}^{\infty}\frac{\sigma_{m}\sigma_{j}\left\langle \phi_{m}\phi_{j},\phi_{n}\right\rangle }{\left|\lambda_{m}\right|-\left|\lambda_{m}\right|-\left|\lambda_{j}\right|} \left(e^{-(\left|\lambda_{m}\right|+\left|\lambda_{j}\right|)T}-e^{-\left|\lambda_{n}\right|T}\right).
\end{align*}
 Note, $|\lambda_{n}|-|\lambda_{m}|-|\lambda_{j}|\neq0$ due to the linear independence of the orthonormal eigenfunctions $\phi_{n}(X)$.  If we remove the terms of order $Te^{-|\lambda_{n}|T}$ as $T\rightarrow\infty$, then we have that 
\[
\dot{\sigma}_{n}=(1-2\sigma_{0}\phi_{0})\sigma_{n},
\]
 or 
\[
\sigma_{n}(\tau)=\frac{\sigma_{n}(0)e^{t}}{\left(1+\sigma_{0}(0)\phi_{0}(e^{t}-1)\right)^{2}}.
\]

Thus, on timescales $T=\mathcal{O}(1/\epsilon)$, or $t=\mathcal{O}(1)$, one has 
\[
u(x,T,t)=u_{0}(x,T,t)+\mathcal{O}(\epsilon).
\]
 From this, since for $n\geq1$, one has by orthonormality 
\[
\int_{-L}^{L}\phi_{n}(x)dx=0,
\]
 it follows that the average of $u(x,T,t)$ to leading order is given by $\sigma_{0}(\tau)\phi_{0}$, or 
\begin{equation}
\left<u\right>\sim\sigma_{0}(t)\phi_{0}=\frac{1}{1+\frac{L-x_{c}}{L+x_{c}}e^{-t}}.\label{avgresult}
\end{equation}
 The quantity $\left<u\right>$ is given by 
\[
\left<u(\cdot,T,t)\right>=\frac{1}{2L}\int_{-L}^{L}u(\xi,T,t)d\xi.
\]
This result shows that the average should increase exponentially fast on time scales of $\mathcal{O}(1)$ or less, which amounts to showing that the transition region of the solution $u$ becomes smoother and ``softer" rapidly.  Having done this, the solution then enters the SFA regime.  
\section{\label{sec:BeyondBarrier}Front Beyond the Turning Point: reduction to a stationary equation}
\label{Reduction to a Stationary Equation}

Since the SFA only leads to a model that describes propagation of a front up to the turning point, we must now find some other means of trying to describe propagation of the front past the turning point.  To do this, by taking $x\gg 1$, so that $a(x)\sim x^{2}$, we suppose that $u=u(\eta(x,t))$, which means the FK equation becomes
\[
\eta_{t}\frac{du}{d\eta}= (x\eta_{x})^{2}\frac{d^{2}u}{d\eta^{2}} + (x^{2}\eta_{x})_{x} \frac{du}{d\eta} + f(u).
\]
By choosing 
\[
\eta_{x} = \pm \frac{1}{x}, ~ \eta_{t} = c, 
\]
we get the stationary FK equation 
\[
\frac{d^{2}u}{d\eta^{2}} + (-c \pm1)  \frac{du}{d\eta} + f(u) = 0,
\]
with coordinate
\[
\eta(x,t) = \pm \ln|x| + ct~.
\]
The variable $\eta$ is another instance of a (TWC), and so again we see that strongly varying diffusion requires a generalization of the definition of a traveling front.  Linearizing around both the ``$+$" and ``$-$" cases from above, for the ``$+$" case, we get roots to the characteristic equation of the form
\[
\lambda^{+}(c) = \frac{(c-1) \pm ((c-1)^{2} - 4)^{1/2}}{2},
\]
and for the ``$-$" case, we get
\[
\lambda^{-}(c) = \frac{(c+1)\pm((c+1)^{2}-4)^{1/2}}{2}.
\]
The front behaves like $u \sim \eta e^{\lambda \eta}$, when $\lambda$ is a double root or $u \sim e^{\lambda \eta}$, when $\lambda$ is a single root.  In either case, we see that the decay rate of the front is not exponentially fast, but is instead only algebraically fast since
\[
e^{\lambda \eta} = |x|^{\pm \lambda} e^{\lambda ct}.  
\]
In order to ensure non-oscillatory decay, in the ``$+$" case, one must take $c<-1$, and in the ``$-$" case, one must take $c>1$.  We define the characteristic curves of the TWC as the curves of constant $\eta(x,t)$ which are given by  
\[
\pm \ln|x_{0}| = \pm\ln|x| + ct, 
\]
so that $x(t) = x_{0}e^{\mp ct}$.  Since the characteristics are identical in both of the relevant cases, we take the ``$-$" case as our convention.  We also see that, in contrast to the characteristic curves in the soft front case, that information is transported away from the origin.  Thus, it appears the TWC can provide a mechanism for propagation of fronts past the turning point.    

However, the algebraic decay is unexpected, and at first glance would seem to imply the TWC does not provide relevant information about the propagation of true fronts, {\it i.e.}, it is not clear that 
\begin{equation}
\lim_{t \rightarrow \infty} u(x,t) = \tilde{u}(\eta(x,t)),
\label{asympconv}
\end{equation}
where $u(x,t)$ is some solution to the FK equation with arbitrary initial condition and $\tilde{u}(\eta)$ denotes a solution to the stationary FK equation.  It is non-trivial, even for the constant diffusion case, to show that solutions of the FK equation for some class of initial data satisfy \eqref{asympconv}.  The literature on this issue is large, so we only refer the reader to the foundational papers \cite{aron, fife, kolm, mckean} as an introduction to this issue. Lacking any rigorous proof of \eqref{asympconv}, we provide a formal result to support the hypothesis that fronts evolve according to the TWC argument.

Motivated then by classical approaches to the FK equation (\cf \cite{murray}), we again look at the linearized FK equation \eqref{linear}.  Taking $x\gg 1$, so that $a(x)\sim x^{2}$, separation of variables can be used to solve \eqref{linear}.  This gives the solution
\[
v(x,t) \sim x^{-\frac{1}{2}\pm \frac{\sqrt{4\bar{c}-3}}{2}} e^{\bar{c}t} = e^{\left(-\frac{1}{2}\pm \frac{\sqrt{4\bar{c}-3}}{2} \right)\ln(x) + \bar{c}t}, 
\]
where we take $\bar{c}\geq 3/4$ so as to eliminate oscillatory solutions.  Thus, we see algebraically decaying fronts are quite natural in this problem. This analysis supports the argument that solutions found via the TWC are in fact representative of true front behavior. Again, this is to say that \eqref{asympconv} holds for some class of initial data.  We also note the TWC fits into the framework of the WKB analysis presented in Section \ref{softfront} since the characteristics past the origin found via the WKB analysis are the same as the TWC with $\eta$ held constant.  

\section{\label{Numerical Results} Numerical Results}

We have developed an explanation for how a strongly varying diffusion coefficient affects the propagation of a front up to the turning point of $a(x)$.  Using the TWC, we have some analytical results explaining front propagation past the turning point.  In this section, we corroborate our results with numerical simulations.  To simulate solutions to the FK equation, we implemented the numerical scheme presented in \cite{branco}, an implicit/explicit method that is second order in space and first order in time.  In these simulations we chose the nonlinearity $f(u) = u(1-u)$, as well as Neumann boundary conditions $u_{x}(\pm L, t) = 0$, and step initial condition as in \eqref{step}.  

\begin{figure}[!h]
\centering
\includegraphics[width=.35\textwidth]{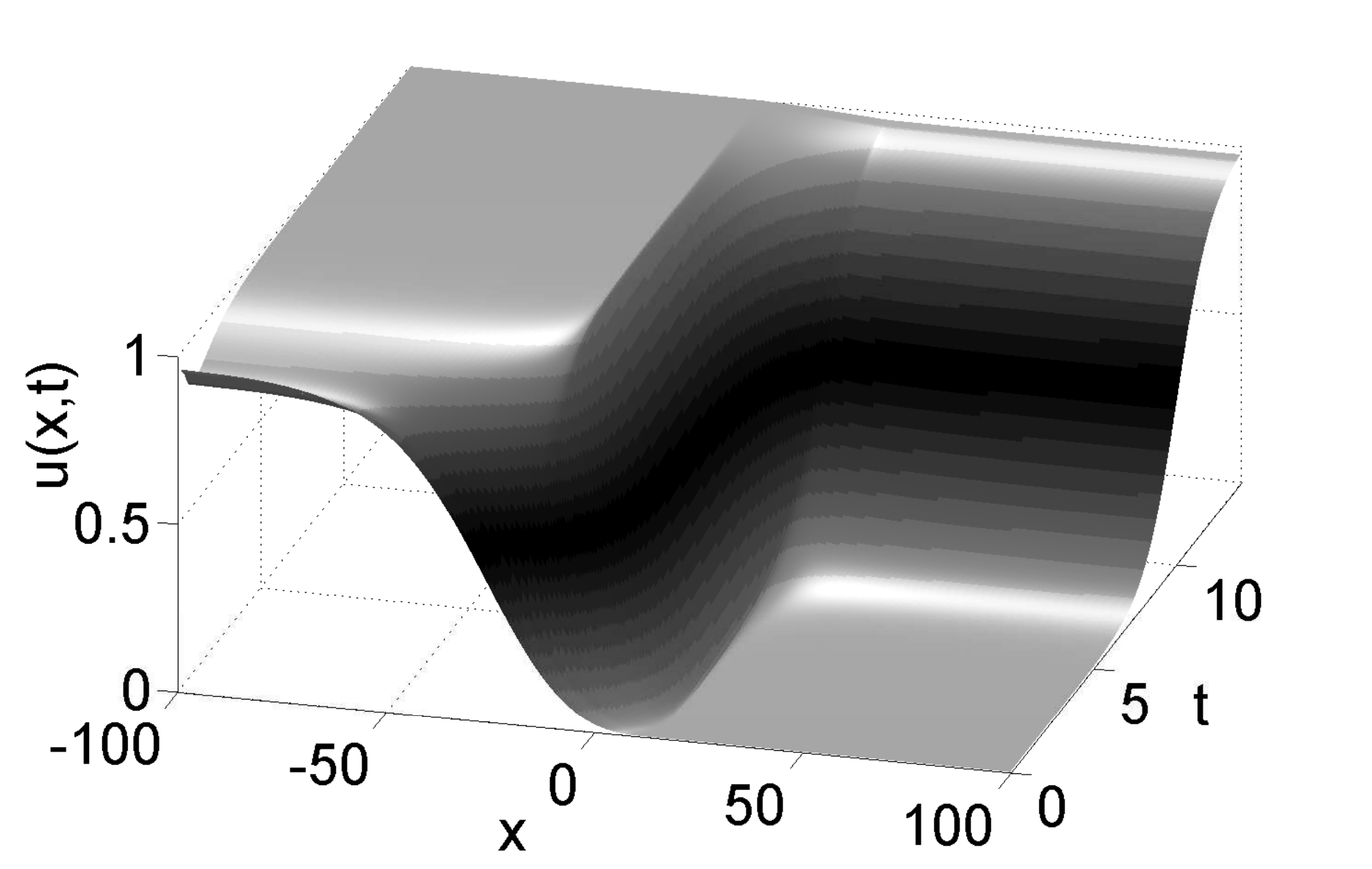}
\caption{ Plot of Solution $u(x,t)$ - $L=100$, $x_{c}(0)=-35$, $\epsilon = .1$}
\label{fig:front_prop}
\end{figure}

\begin{figure}[!h]
\centering
\includegraphics[width=.35\textwidth]{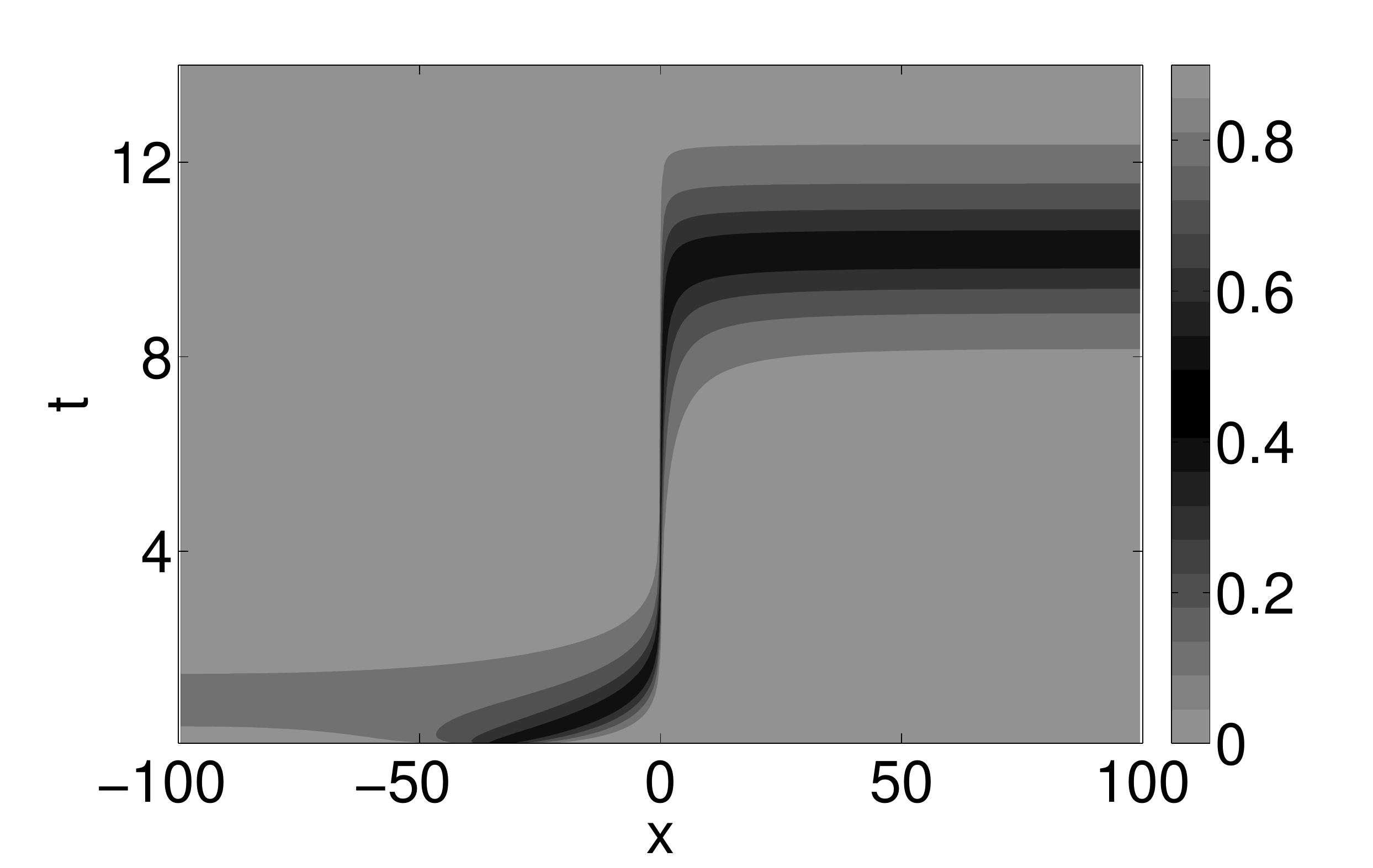}
\caption{Contour Plot of Solution $u(x,t)$ - $L=100$, $x_{c}(0)=-35$, $\epsilon = .1$}
\label{fig:front_prop_bird}
\end{figure}

\begin{figure}[!h]
\centering
\includegraphics[width=.35\textwidth]{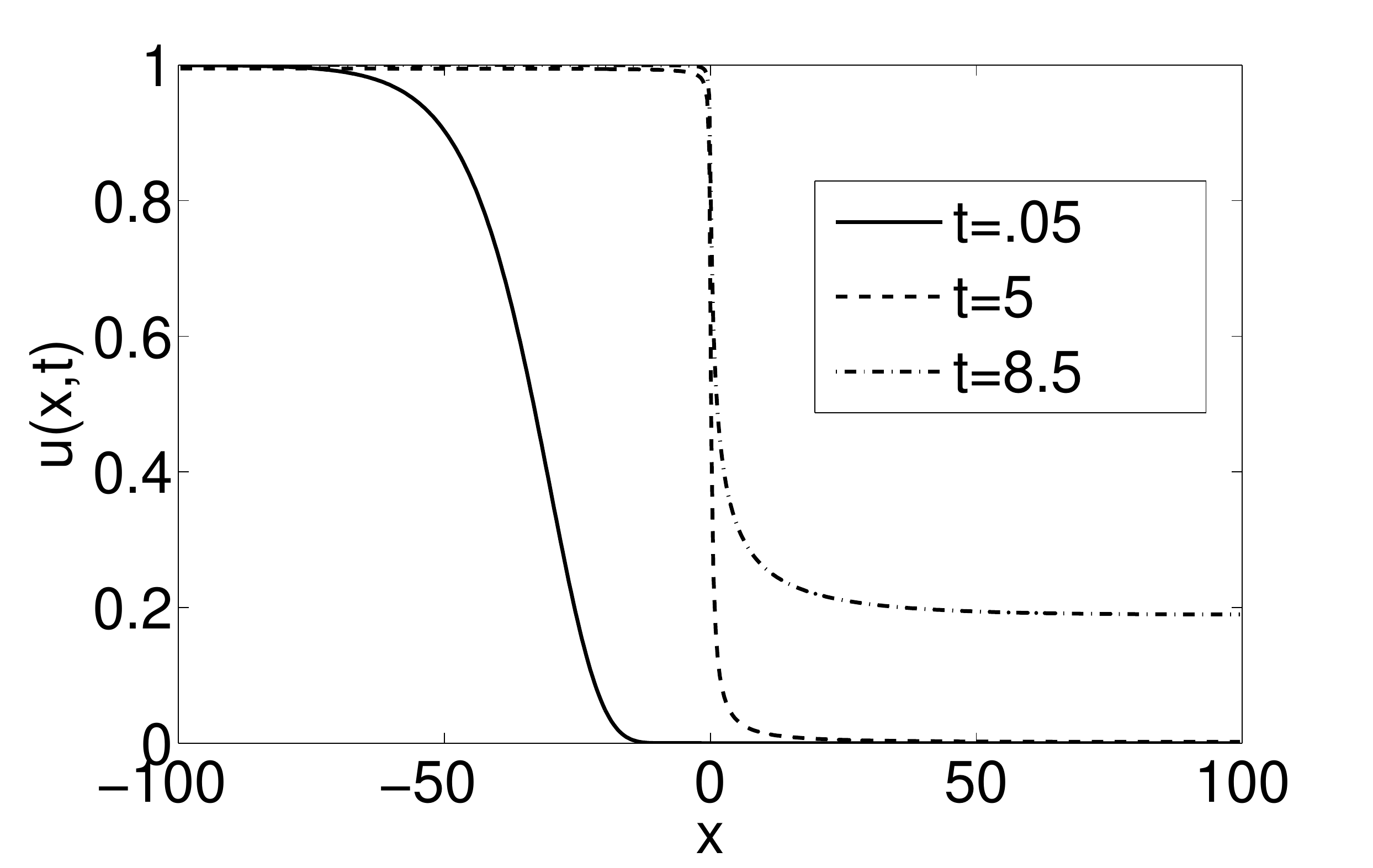}
\caption{Solution Profiles at Various Times - $L=100$, $x_{c}(0)=-35$, $\epsilon = .1$}
\label{fig:front_prop_profile}
\end{figure}

\begin{figure}[!h]
\centering
\includegraphics[width=.35\textwidth]{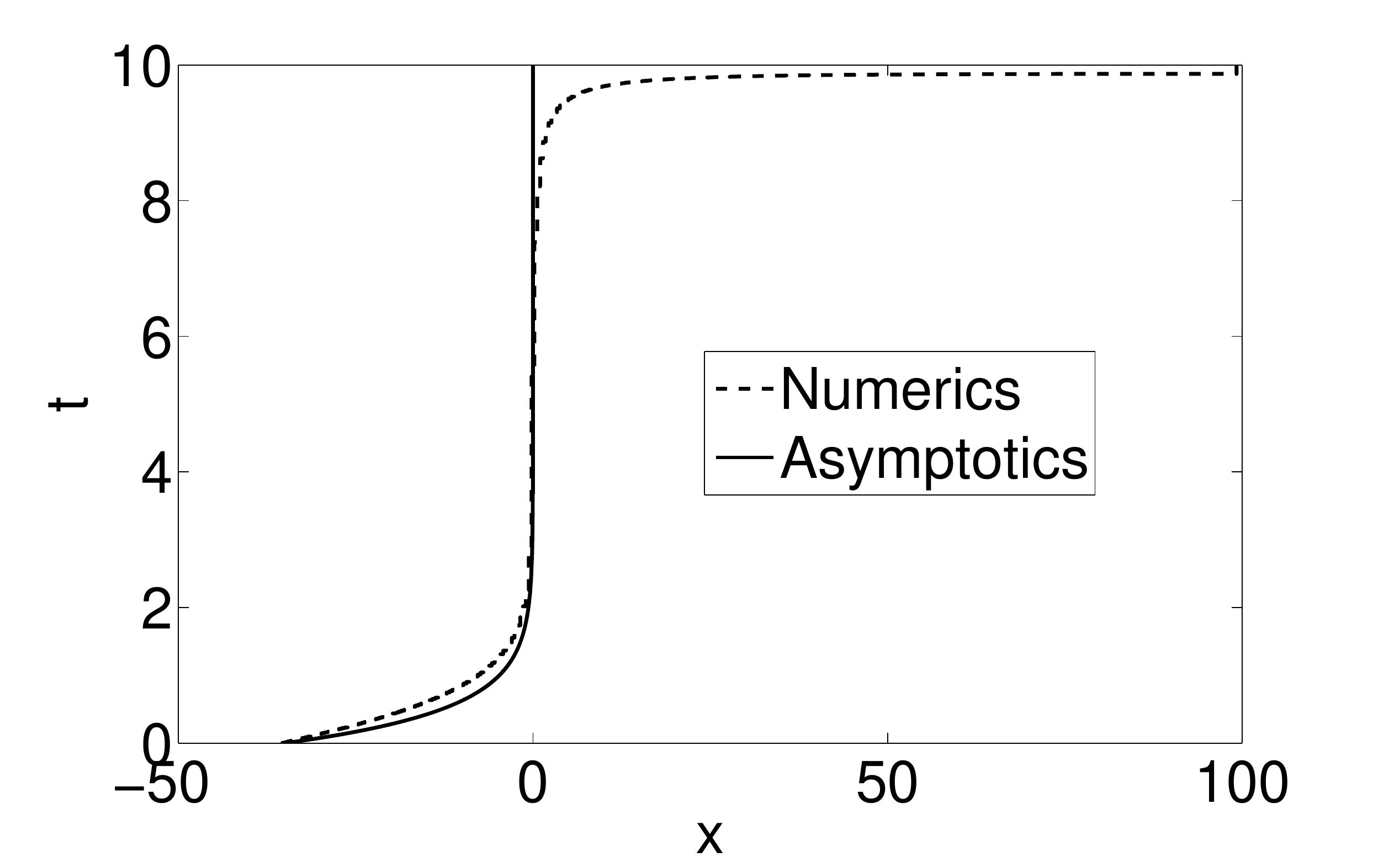}
\caption{Comparison of SFA to Numerics - $L=100$, $x_{c}(0)=-35$, $\epsilon = .1$}
\label{fig:front_track}
\end{figure}

Figures \ref{fig:front_prop} and \ref{fig:front_prop_bird}
depict the development of the solution in the $(t,x)$ plane as a surface and as a contour plot, respectively.  Figure \ref{fig:front_prop_profile} 
shows the front before, at, and after the turning point.  In Figure \ref{fig:front_track}, 
we compare the numerical propagation of the front and our asymptotic theory.  The two curves represent the function $x_{c,n}(t)$ at which $u_{num}(x_{c,n}(t),t)=1/2$, where $u_{num}(x,t)$ denotes the numerical approximation to $u(x,t)$.   In the case of the asymptotic curve, we take $t_{0}=0$, and use the same step initial condition as used in the numerics in \eqref{sfasol}.  The agreement between the SFA and numerics is convincing, and thus we have confirmation that the SFA is valid up to the turning point.  As can be seen the SFA breaks down at the turning point.  From Figures \ref{fig:front_prop}, \ref{fig:front_prop_bird}, and \ref{fig:front_prop_profile}, we see that the front steepens substantially.    

\begin{figure}[!h]
\centering
\includegraphics[width=.35\textwidth]{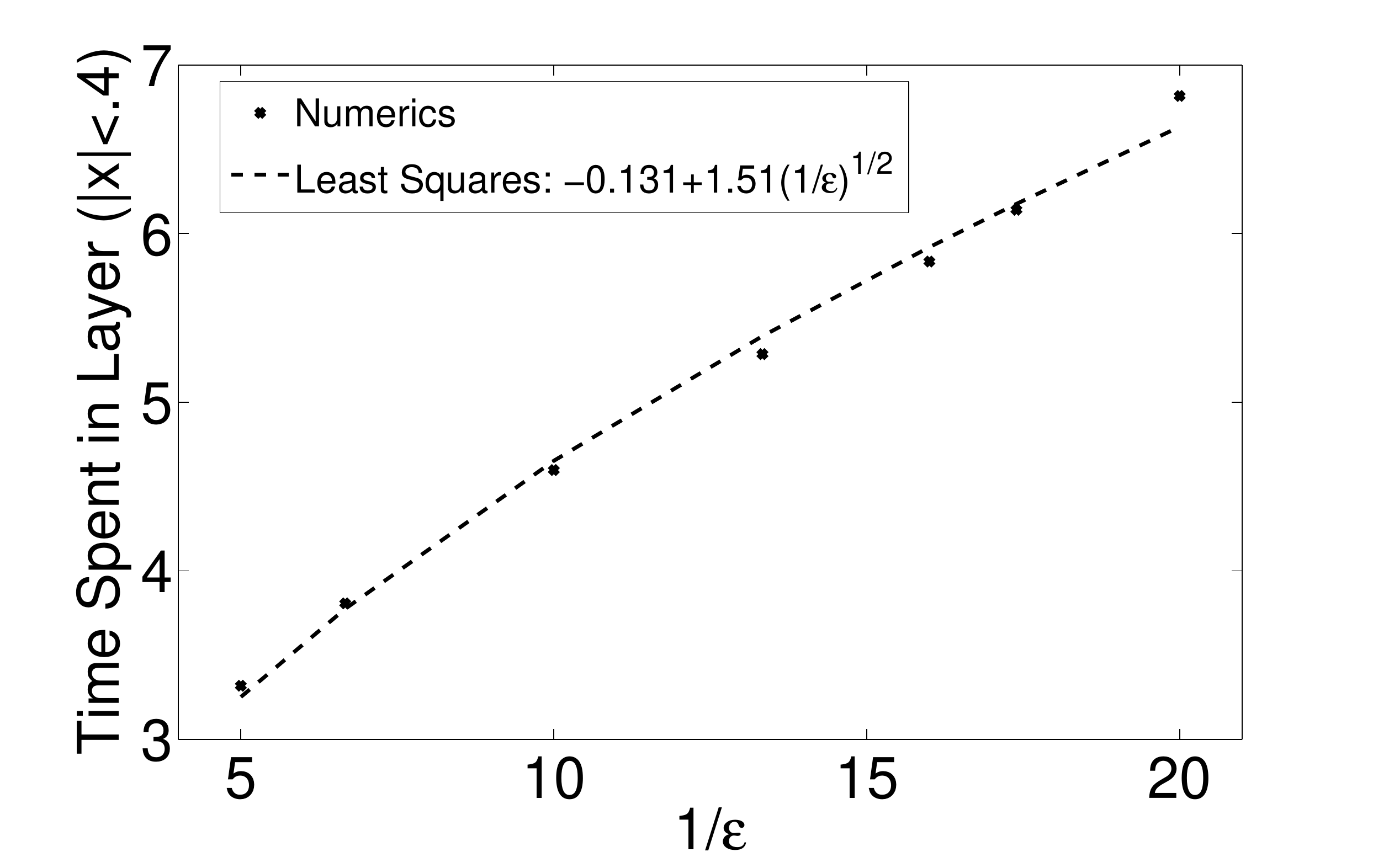}
\caption{Trapping Times in Layer Around Minimum - $L=100$, $x_{c}(0)=-35$}
\label{fig:trap_times}
\end{figure}

Once the front reaches the turning point, we see in Figures \ref{fig:front_prop}, \ref{fig:front_prop_bird} and \ref{fig:front_track} that the front becomes trapped.  In Section \ref{softfront}, we have predicted the timescale of trapping to be $\mathcal{O}(1/\sqrt{\epsilon})$.  In Figure \ref{fig:trap_times}, 
we plot the time duration that the simulated front spends near $x=0$ as a function of $1/\epsilon$, where we define {\it near} as $|x|<0.4$.  Note the value $0.4$ is the width of the spatial mesh step in the numerics, and thus the smallest scale on which phenomena can be distinguished.  Figure \ref{fig:trap_times} shows that a least squares fit of the times to a curve growing like $1/\sqrt{\epsilon}$ is accurate.  Note that a decrease the spatial mesh step size will reduce the error between observed time and the least squares fit curve (results not presented here).

On the other side of the turning point, as seen in Figure \ref{fig:front_prop_bird}, the solution propagates in a manner consistent with the log-transformed TWC, {\it i.e.}, information is propagating along curves of the form $\eta = -\ln|x| + ct$, or $x = x_{0}e^{ct}$.  However, it is currently an open question as to how $c$ is selected (as opposed to the constant coefficient case where the regularity of the initial condition chooses $c$) and this is a direction of future work.
\section{Conclusions and Future Work}

We have shown in the paper that, after diffusion smooths the initial conditions such that SFA is valid, the choice of a strongly varying diffusion coefficient with a global minimum implies 
\begin{itemize}
\item The definition of a traveling front must be generalized via the {\it Traveling Wave Coordinate} as defined in Section \ref{sec:SFC}. 
\item Using WKB analysis, we show the minimum, or {\it turning point}, of $a(x)$ causes the formation of shock-like behavior and leads to a trapping of the front on asymptotically long time scales. 
\item The behavior of the front on either side of the turning point is fundamentally different, and on either side, the TWC is necessary to describe dynamics.  
\item After the front is past the turning point, the TWC allows us to transform into a stationary FK equation, thus simplifying all subsequent analysis.  
\end{itemize}

To the best of our knowledge, the case of strongly varying diffusion has not been studied in the literature, and this article is the first investigation of traveling waves with continuously and strongly varying spatial diffusion.  We think these results will prove especially useful in modeling real world data for two reasons.

First, quadratic spatial diffusion could be used to model barriers such as mountains and our results will give estimates of a time to barrier transit.  For example, the invasion of the midwestern US by gypsy moths has been extensively studied (\cf \cite{Tobin2007} and references therein).  It is also well known that the moths' spatial diffusion rate is tightly correlated with the local habitat \cite{Liebhold2008}.  By estimating the population growth rates and fitting a quadratic curve to the local diffusion rates from capture-mark-recapture experiments (like those described in \cite{Andow1990}), the method developed in this paper could be used to estimate the waiting time until the moths reach a certain habitat.  In a non-ecological context such as the morphogen gradients discussed in \cite{Philip1992}, these {\it barriers} could be nutrient deficient regions in gap junctions.  Here, our results could be used to estimate the time till the concentration of morphogen or metabolite at some boundary reaches a critical threshold.

Second, with our work, we now have a non-trivial spatially varying FK equation that (via a traveling wave ansatz) could be reduced to a much simpler stationary problem.  This will allow for a greater degree of control and flexibility in curve fitting routines as it will only require a snapshot at a single point in time.  With the gyspsy moth invasion example, one could envision estimating the wave speed from a single snapshot.  Furthermore, the authors currently have a manuscript in progress applying these tools to the gypsy moth invasion in an effort to estimate the time to barrier transit and invasion wave speed.

Naturally, there are also many future directions for the analytical aspect of this work.  It would be of particular interest to develop asymptotic matching techniques that would allow for connecting the various regimes of the front, {\it i.e.}, the sharp to soft front transition, and then the trapping layer around the minimum of $a(x)$.  Likewise, it would also be of interest to develop higher order expansions for the speed $c$ in the SFA regime.     
\bibliography{latest_fisher}

\begin{thebibliography}{40}
\expandafter\ifx\csname natexlab\endcsname\relax\def\natexlab#1{#1}\fi
\expandafter\ifx\csname bibnamefont\endcsname\relax
  \def\bibnamefont#1{#1}\fi
\expandafter\ifx\csname bibfnamefont\endcsname\relax
  \def\bibfnamefont#1{#1}\fi
\expandafter\ifx\csname citenamefont\endcsname\relax
  \def\citenamefont#1{#1}\fi
\expandafter\ifx\csname url\endcsname\relax
  \def\url#1{\texttt{#1}}\fi
\expandafter\ifx\csname urlprefix\endcsname\relax\def\urlprefix{URL }\fi
\providecommand{\bibinfo}[2]{#2}
\providecommand{\eprint}[2][]{\url{#2}}

\bibitem[{\citenamefont{Murray}(2001)}]{murray}
\bibinfo{author}{\bibfnamefont{J.~D.} \bibnamefont{Murray}},
  \emph{\bibinfo{title}{{Mathematical Biology I: An Introduction}}},
  vol.~\bibinfo{volume}{17} of \emph{\bibinfo{series}{Interdisciplinary Applied
  Math.}} (\bibinfo{publisher}{Springer-Verlag}, \bibinfo{address}{New York,
  NY}, \bibinfo{year}{2001}), \bibinfo{edition}{3rd} ed.

\bibitem[{\citenamefont{Murray}(2003)}]{Murray2003}
\bibinfo{author}{\bibfnamefont{J.~D.} \bibnamefont{Murray}},
  \emph{\bibinfo{title}{{Mathematical Biology II: Spatial Models and Biomedical
  Applications}}}, vol.~\bibinfo{volume}{18} of
  \emph{\bibinfo{series}{Interdisciplinary Applied Math.}}
  (\bibinfo{publisher}{Springer-Verlag}, \bibinfo{address}{New York, NY},
  \bibinfo{year}{2003}), \bibinfo{edition}{3rd} ed.

\bibitem[{\citenamefont{Shigesada and Kawasaki}(1997)}]{shigesada}
\bibinfo{author}{\bibfnamefont{N.}~\bibnamefont{Shigesada}} \bibnamefont{and}
  \bibinfo{author}{\bibfnamefont{K.}~\bibnamefont{Kawasaki}},
  \emph{\bibinfo{title}{Biological Invasions: Theory and Practice}}
  (\bibinfo{publisher}{Oxford University Press}, \bibinfo{address}{New York,
  NY}, \bibinfo{year}{1997}).

\bibitem[{\citenamefont{Maini et~al.}(2004)\citenamefont{Maini, Mcelwain, and
  Leavesley}}]{Maini2004}
\bibinfo{author}{\bibfnamefont{P.~K.} \bibnamefont{Maini}},
  \bibinfo{author}{\bibfnamefont{D.~L.~S.} \bibnamefont{Mcelwain}},
  \bibnamefont{and}
  \bibinfo{author}{\bibfnamefont{D.}~\bibnamefont{Leavesley}},
  \bibinfo{journal}{Applied Mathematics Letters} \textbf{\bibinfo{volume}{17}},
  \bibinfo{pages}{575} (\bibinfo{year}{2004}).

\bibitem[{\citenamefont{Epstein and Pojman}(1998)}]{epstein}
\bibinfo{author}{\bibfnamefont{I.}~\bibnamefont{Epstein}} \bibnamefont{and}
  \bibinfo{author}{\bibfnamefont{J.~A.} \bibnamefont{Pojman}},
  \emph{\bibinfo{title}{An Introduction to Nonlinear Chemical Dynamics}}
  (\bibinfo{publisher}{Oxford University Press}, \bibinfo{address}{New York,
  N.Y.}, \bibinfo{year}{1998}).

\bibitem[{\citenamefont{Grindrod}(1996)}]{grindrod}
\bibinfo{author}{\bibfnamefont{P.}~\bibnamefont{Grindrod}},
  \emph{\bibinfo{title}{The Theory and Applications of Reaction-Diffusion
  Equations}} (\bibinfo{publisher}{Oxford University Press},
  \bibinfo{address}{New York, NY}, \bibinfo{year}{1996}).

\bibitem[{\citenamefont{Hammond and Bortz}(2011)}]{HammondBortz2011}
\bibinfo{author}{\bibfnamefont{J.~F.} \bibnamefont{Hammond}} \bibnamefont{and}
  \bibinfo{author}{\bibfnamefont{D.~M.} \bibnamefont{Bortz}},
  \bibinfo{journal}{Applied Mathematics and Computation}
  \textbf{\bibinfo{volume}{218}}, \bibinfo{pages}{2497} (\bibinfo{year}{2011}).

\bibitem[{\citenamefont{Shebesch and Engel}(1998)}]{engel}
\bibinfo{author}{\bibfnamefont{I.}~\bibnamefont{Shebesch}} \bibnamefont{and}
  \bibinfo{author}{\bibfnamefont{H.}~\bibnamefont{Engel}},
  \bibinfo{journal}{Phys. Rev. E.} \textbf{\bibinfo{volume}{57}},
  \bibinfo{pages}{3905} (\bibinfo{year}{1998}).

\bibitem[{\citenamefont{Cuesta and King}(2010)}]{cuesta}
\bibinfo{author}{\bibfnamefont{C.}~\bibnamefont{Cuesta}} \bibnamefont{and}
  \bibinfo{author}{\bibfnamefont{J.}~\bibnamefont{King}},
  \bibinfo{journal}{Q.Jl Mech.Appl.Math} \textbf{\bibinfo{volume}{63}},
  \bibinfo{pages}{521} (\bibinfo{year}{2010}).

\bibitem[{\citenamefont{Keener}(2000)}]{keener}
\bibinfo{author}{\bibfnamefont{J.~P.} \bibnamefont{Keener}},
  \bibinfo{journal}{SIAM J. Appl. Math.} \textbf{\bibinfo{volume}{81}},
  \bibinfo{pages}{317} (\bibinfo{year}{2000}).

\bibitem[{\citenamefont{M\'{e}ndez et~al.}(2003)\citenamefont{M\'{e}ndez, Fort,
  Rotstein, and Fedotov}}]{mendez}
\bibinfo{author}{\bibfnamefont{V.}~\bibnamefont{M\'{e}ndez}},
  \bibinfo{author}{\bibfnamefont{J.}~\bibnamefont{Fort}},
  \bibinfo{author}{\bibfnamefont{H.}~\bibnamefont{Rotstein}}, \bibnamefont{and}
  \bibinfo{author}{\bibfnamefont{S.}~\bibnamefont{Fedotov}},
  \bibinfo{journal}{Phys. Rev. E} \textbf{\bibinfo{volume}{68}},
  \bibinfo{pages}{041105} (\bibinfo{year}{2003}).

\bibitem[{\citenamefont{Xin}(2000)}]{xin}
\bibinfo{author}{\bibfnamefont{J.}~\bibnamefont{Xin}}, \bibinfo{journal}{SIAM
  Review} \textbf{\bibinfo{volume}{42}}, \bibinfo{pages}{161}
  (\bibinfo{year}{2000}).

\bibitem[{\citenamefont{Andow et~al.}(1990)\citenamefont{Andow, Kareiva, Levin,
  and Okubo}}]{Andow1990}
\bibinfo{author}{\bibfnamefont{D.~A.} \bibnamefont{Andow}},
  \bibinfo{author}{\bibfnamefont{P.~M.} \bibnamefont{Kareiva}},
  \bibinfo{author}{\bibfnamefont{S.~A.} \bibnamefont{Levin}}, \bibnamefont{and}
  \bibinfo{author}{\bibfnamefont{A.}~\bibnamefont{Okubo}},
  \bibinfo{journal}{Landscape Ecology} \textbf{\bibinfo{volume}{4}},
  \bibinfo{pages}{177} (\bibinfo{year}{1990}).

\bibitem[{\citenamefont{Cruywagen et~al.}(1996)\citenamefont{Cruywagen,
  Kareiva, Lewis, and Murray}}]{Cruywagen1996}
\bibinfo{author}{\bibfnamefont{G.~C.} \bibnamefont{Cruywagen}},
  \bibinfo{author}{\bibfnamefont{P.}~\bibnamefont{Kareiva}},
  \bibinfo{author}{\bibfnamefont{M.~a.} \bibnamefont{Lewis}}, \bibnamefont{and}
  \bibinfo{author}{\bibfnamefont{J.~D.} \bibnamefont{Murray}},
  \bibinfo{journal}{Theoretical population biology}
  \textbf{\bibinfo{volume}{49}}, \bibinfo{pages}{1} (\bibinfo{year}{1996}).

\bibitem[{\citenamefont{Dobzhansky et~al.}(1979)\citenamefont{Dobzhansky,
  Powell, Taylor, and Andregg}}]{Dobzhansky1979}
\bibinfo{author}{\bibfnamefont{T.}~\bibnamefont{Dobzhansky}},
  \bibinfo{author}{\bibfnamefont{J.~R.} \bibnamefont{Powell}},
  \bibinfo{author}{\bibfnamefont{C.~E.} \bibnamefont{Taylor}},
  \bibnamefont{and} \bibinfo{author}{\bibfnamefont{M.}~\bibnamefont{Andregg}},
  \bibinfo{journal}{The American Naturalist} \textbf{\bibinfo{volume}{114}},
  \bibinfo{pages}{325} (\bibinfo{year}{1979}).

\bibitem[{\citenamefont{Hastings}(1982)}]{Hastings1982}
\bibinfo{author}{\bibfnamefont{A.}~\bibnamefont{Hastings}},
  \bibinfo{journal}{Journal of Mathematical Biology}
  \textbf{\bibinfo{volume}{16}}, \bibinfo{pages}{49} (\bibinfo{year}{1982}).

\bibitem[{\citenamefont{Turchin and Thoeny}(1993)}]{Turchin1993}
\bibinfo{author}{\bibfnamefont{P.}~\bibnamefont{Turchin}} \bibnamefont{and}
  \bibinfo{author}{\bibfnamefont{W.~T.} \bibnamefont{Thoeny}},
  \bibinfo{journal}{Ecological Applications} \textbf{\bibinfo{volume}{3}},
  \bibinfo{pages}{187} (\bibinfo{year}{1993}).

\bibitem[{\citenamefont{Harrison}(1980)}]{Harrison1980}
\bibinfo{author}{\bibfnamefont{R.~G.} \bibnamefont{Harrison}},
  \bibinfo{journal}{Annual Review of Ecology and Systematics}
  \textbf{\bibinfo{volume}{11}}, \bibinfo{pages}{95} (\bibinfo{year}{1980}).

\bibitem[{\citenamefont{Slatkin}(1985)}]{Slatkin1985}
\bibinfo{author}{\bibfnamefont{M.}~\bibnamefont{Slatkin}},
  \bibinfo{journal}{Annual Review of Ecology and Systematics}
  \textbf{\bibinfo{volume}{16}}, \bibinfo{pages}{393} (\bibinfo{year}{1985}).

\bibitem[{\citenamefont{Bridle et~al.}(2001)\citenamefont{Bridle, Baird, and
  Butlin}}]{Bridle2001}
\bibinfo{author}{\bibfnamefont{J.~R.} \bibnamefont{Bridle}},
  \bibinfo{author}{\bibfnamefont{S.~J.} \bibnamefont{Baird}}, \bibnamefont{and}
  \bibinfo{author}{\bibfnamefont{R.~K.} \bibnamefont{Butlin}},
  \bibinfo{journal}{Evolution} \textbf{\bibinfo{volume}{55}},
  \bibinfo{pages}{1832} (\bibinfo{year}{2001}).

\bibitem[{\citenamefont{Kinezaki et~al.}(2003)\citenamefont{Kinezaki, Kawasaki,
  Takasu, and Shigesada}}]{Kinezaki2003}
\bibinfo{author}{\bibfnamefont{N.}~\bibnamefont{Kinezaki}},
  \bibinfo{author}{\bibfnamefont{K.}~\bibnamefont{Kawasaki}},
  \bibinfo{author}{\bibfnamefont{F.}~\bibnamefont{Takasu}}, \bibnamefont{and}
  \bibinfo{author}{\bibfnamefont{N.}~\bibnamefont{Shigesada}},
  \bibinfo{journal}{Theoretical Population Biology}
  \textbf{\bibinfo{volume}{64}}, \bibinfo{pages}{291} (\bibinfo{year}{2003}).

\bibitem[{\citenamefont{Kinezaki et~al.}(2010)\citenamefont{Kinezaki, Kawasaki,
  and Shigesada}}]{Kinezaki2010}
\bibinfo{author}{\bibfnamefont{N.}~\bibnamefont{Kinezaki}},
  \bibinfo{author}{\bibfnamefont{K.}~\bibnamefont{Kawasaki}}, \bibnamefont{and}
  \bibinfo{author}{\bibfnamefont{N.}~\bibnamefont{Shigesada}},
  \bibinfo{journal}{Theoretical population biology}
  \textbf{\bibinfo{volume}{78}}, \bibinfo{pages}{298} (\bibinfo{year}{2010}).

\bibitem[{\citenamefont{Shigesada et~al.}(1986)\citenamefont{Shigesada,
  Kawasaki, and Teramoto}}]{Shigesada1986}
\bibinfo{author}{\bibfnamefont{N.}~\bibnamefont{Shigesada}},
  \bibinfo{author}{\bibfnamefont{K.}~\bibnamefont{Kawasaki}}, \bibnamefont{and}
  \bibinfo{author}{\bibfnamefont{E.}~\bibnamefont{Teramoto}},
  \bibinfo{journal}{Theoretical Population Biology}
  \textbf{\bibinfo{volume}{30}}, \bibinfo{pages}{143} (\bibinfo{year}{1986}).

\bibitem[{\citenamefont{Petrovskii and Li}(2005)}]{Petrovskii2005}
\bibinfo{author}{\bibfnamefont{S.}~\bibnamefont{Petrovskii}} \bibnamefont{and}
  \bibinfo{author}{\bibfnamefont{B.-L.} \bibnamefont{Li}},
  \emph{\bibinfo{title}{{Exactly Solvable Models of Biological Invasion}}},
  vol.~\bibinfo{volume}{7} of \emph{\bibinfo{series}{Chapman and Hall/CRC
  Mathematical and Computational Biology}} (\bibinfo{publisher}{Chapman and
  Hall/CRC}, \bibinfo{year}{2005}).

\bibitem[{\citenamefont{Othmer and Pate}(1980)}]{Othmer1980}
\bibinfo{author}{\bibfnamefont{H.~G.} \bibnamefont{Othmer}} \bibnamefont{and}
  \bibinfo{author}{\bibfnamefont{E.}~\bibnamefont{Pate}},
  \bibinfo{journal}{Proceedings of the National Academy of Sciences of the
  United States of America} \textbf{\bibinfo{volume}{77}},
  \bibinfo{pages}{4180} (\bibinfo{year}{1980}).

\bibitem[{\citenamefont{Maini et~al.}(1992)\citenamefont{Maini, Benson, and
  Sherratt}}]{Philip1992}
\bibinfo{author}{\bibfnamefont{P.~K.} \bibnamefont{Maini}},
  \bibinfo{author}{\bibfnamefont{D.~L.} \bibnamefont{Benson}},
  \bibnamefont{and} \bibinfo{author}{\bibfnamefont{J.~A.}
  \bibnamefont{Sherratt}}, \bibinfo{journal}{Mathematical Medicine and Biology}
  \textbf{\bibinfo{volume}{9}}, \bibinfo{pages}{197} (\bibinfo{year}{1992}).

\bibitem[{\citenamefont{Benson et~al.}(1993)\citenamefont{Benson, Maini, and
  Sherratt}}]{Benson1993}
\bibinfo{author}{\bibfnamefont{D.}~\bibnamefont{Benson}},
  \bibinfo{author}{\bibfnamefont{P.}~\bibnamefont{Maini}}, \bibnamefont{and}
  \bibinfo{author}{\bibfnamefont{J.}~\bibnamefont{Sherratt}}, in
  \emph{\bibinfo{booktitle}{Experimental and Theoretical Advances in Biological
  Pattern Formation}}, edited by \bibinfo{editor}{\bibfnamefont{H.~G.}
  \bibnamefont{Othmer}}, \bibinfo{editor}{\bibfnamefont{P.~K.}
  \bibnamefont{Maini}}, \bibnamefont{and} \bibinfo{editor}{\bibfnamefont{J.~D.}
  \bibnamefont{Murray}} (\bibinfo{publisher}{Plenum Pres},
  \bibinfo{year}{1993}), Nato Science Series: A:, pp. \bibinfo{pages}{29--32}.

\bibitem[{\citenamefont{Chaplain et~al.}(1994)\citenamefont{Chaplain, Benson,
  and Maini}}]{Chaplain1994}
\bibinfo{author}{\bibfnamefont{M.~A.} \bibnamefont{Chaplain}},
  \bibinfo{author}{\bibfnamefont{D.~L.} \bibnamefont{Benson}},
  \bibnamefont{and} \bibinfo{author}{\bibfnamefont{P.~K.} \bibnamefont{Maini}},
  \bibinfo{journal}{Mathematical biosciences} \textbf{\bibinfo{volume}{121}},
  \bibinfo{pages}{1} (\bibinfo{year}{1994}).

\bibitem[{\citenamefont{Benson et~al.}(1998)\citenamefont{Benson, Maini, and
  Sherratt}}]{Benson1998}
\bibinfo{author}{\bibfnamefont{D.~L.} \bibnamefont{Benson}},
  \bibinfo{author}{\bibfnamefont{P.~K.} \bibnamefont{Maini}}, \bibnamefont{and}
  \bibinfo{author}{\bibfnamefont{J.~a.} \bibnamefont{Sherratt}},
  \bibinfo{journal}{Journal of Mathematical Biology}
  \textbf{\bibinfo{volume}{37}}, \bibinfo{pages}{381} (\bibinfo{year}{1998}).

\bibitem[{\citenamefont{Aronson and Weinberger}(1975)}]{aron}
\bibinfo{author}{\bibfnamefont{D.}~\bibnamefont{Aronson}} \bibnamefont{and}
  \bibinfo{author}{\bibfnamefont{H.}~\bibnamefont{Weinberger}},
  \bibinfo{journal}{Lecture Notes in Math.} \textbf{\bibinfo{volume}{446}},
  \bibinfo{pages}{5} (\bibinfo{year}{1975}).

\bibitem[{\citenamefont{Booty et~al.}(1993)\citenamefont{Booty, Haberman, and
  Minzoni}}]{haberman}
\bibinfo{author}{\bibfnamefont{M.}~\bibnamefont{Booty}},
  \bibinfo{author}{\bibfnamefont{R.}~\bibnamefont{Haberman}}, \bibnamefont{and}
  \bibinfo{author}{\bibfnamefont{A.}~\bibnamefont{Minzoni}},
  \bibinfo{journal}{SIAM J. Appl. Math.} \textbf{\bibinfo{volume}{53}},
  \bibinfo{pages}{1009} (\bibinfo{year}{1993}).

\bibitem[{\citenamefont{Evans}(1999)}]{evans}
\bibinfo{author}{\bibfnamefont{L.}~\bibnamefont{Evans}},
  \emph{\bibinfo{title}{Partial Differential Equations}}
  (\bibinfo{publisher}{AMS}, \bibinfo{address}{Providence, R.I.},
  \bibinfo{year}{1999}).

\bibitem[{\citenamefont{Guenther and Lee}(1996)}]{guenther}
\bibinfo{author}{\bibfnamefont{R.}~\bibnamefont{Guenther}} \bibnamefont{and}
  \bibinfo{author}{\bibfnamefont{J.}~\bibnamefont{Lee}},
  \emph{\bibinfo{title}{Partial Differential Equations of Mathematical Physics
  and Integral Equations}} (\bibinfo{publisher}{Dover},
  \bibinfo{address}{Mineola, NY}, \bibinfo{year}{1996}).

\bibitem[{\citenamefont{Bender and Orszag}(1999)}]{Bender}
\bibinfo{author}{\bibfnamefont{C.~M.} \bibnamefont{Bender}} \bibnamefont{and}
  \bibinfo{author}{\bibfnamefont{S.}~\bibnamefont{Orszag}},
  \emph{\bibinfo{title}{Advanced Mathematical Methods Methods for Scientists
  and Engineers: Asymptotic Methods and Perturbation Theory}}
  (\bibinfo{publisher}{Springer}, \bibinfo{address}{New York, NY},
  \bibinfo{year}{1999}).

\bibitem[{\citenamefont{Fife}(1979)}]{fife}
\bibinfo{author}{\bibfnamefont{P.}~\bibnamefont{Fife}},
  \emph{\bibinfo{title}{Mathematical Aspects of Reacting and Diffusing
  Systems}} (\bibinfo{publisher}{Springer-Verlag}, \bibinfo{address}{New York,
  NY}, \bibinfo{year}{1979}).

\bibitem[{\citenamefont{Kolmogorov et~al.}(1937)\citenamefont{Kolmogorov,
  Petrovskii, and Piskunov}}]{kolm}
\bibinfo{author}{\bibfnamefont{A.}~\bibnamefont{Kolmogorov}},
  \bibinfo{author}{\bibfnamefont{I.}~\bibnamefont{Petrovskii}},
  \bibnamefont{and} \bibinfo{author}{\bibfnamefont{N.}~\bibnamefont{Piskunov}},
  \bibinfo{journal}{Mosc. Univ. Bull. Math.} \textbf{\bibinfo{volume}{1}},
  \bibinfo{pages}{1} (\bibinfo{year}{1937}).

\bibitem[{\citenamefont{McKean}(1976)}]{mckean}
\bibinfo{author}{\bibfnamefont{H.}~\bibnamefont{McKean}},
  \bibinfo{journal}{Commun. Pure Appl. Math.} \textbf{\bibinfo{volume}{29}},
  \bibinfo{pages}{323} (\bibinfo{year}{1976}).

\bibitem[{\citenamefont{Branco et~al.}(2007)\citenamefont{Branco, Ferreira, and
  de~Oliveira}}]{branco}
\bibinfo{author}{\bibfnamefont{J.}~\bibnamefont{Branco}},
  \bibinfo{author}{\bibfnamefont{J.}~\bibnamefont{Ferreira}}, \bibnamefont{and}
  \bibinfo{author}{\bibfnamefont{P.}~\bibnamefont{de~Oliveira}},
  \bibinfo{journal}{Applied Numerical Mathematics}
  \textbf{\bibinfo{volume}{57}}, \bibinfo{pages}{89} (\bibinfo{year}{2007}).

\bibitem[{\citenamefont{Tobin et~al.}(2007)\citenamefont{Tobin, Whitmire,
  Johnson, Bj\o~rnstad, and Liebhold}}]{Tobin2007}
\bibinfo{author}{\bibfnamefont{P.~C.} \bibnamefont{Tobin}},
  \bibinfo{author}{\bibfnamefont{S.~L.} \bibnamefont{Whitmire}},
  \bibinfo{author}{\bibfnamefont{D.~M.} \bibnamefont{Johnson}},
  \bibinfo{author}{\bibfnamefont{O.~N.} \bibnamefont{Bj\o~rnstad}},
  \bibnamefont{and} \bibinfo{author}{\bibfnamefont{A.~M.}
  \bibnamefont{Liebhold}}, \bibinfo{journal}{Ecology letters}
  \textbf{\bibinfo{volume}{10}}, \bibinfo{pages}{36} (\bibinfo{year}{2007}).

\bibitem[{\citenamefont{Liebhold and Tobin}(2008)}]{Liebhold2008}
\bibinfo{author}{\bibfnamefont{A.~M.} \bibnamefont{Liebhold}} \bibnamefont{and}
  \bibinfo{author}{\bibfnamefont{P.~C.} \bibnamefont{Tobin}},
  \bibinfo{journal}{Annual review of entomology} \textbf{\bibinfo{volume}{53}},
  \bibinfo{pages}{387} (\bibinfo{year}{2008}).

\end{thebibliography}
\bibliographystyle{apsrev}
\end{document}